\documentclass[12pt]{amsart}
\usepackage{amsmath,amscd,amssymb,amsfonts}
\newcommand{\ver}{June 15, 2006, v.2}
\newcommand{\scirc}{\raise.2ex\hbox{${\scriptstyle\circ}$}}
\newcommand{\ssbull}{\raise.2ex\hbox{${\scriptscriptstyle\bullet}$}}
\newcommand{\mopls}{\hbox{$\bigoplus$}}
\newcommand{\bC}{{\mathbb C}}
\newcommand{\bD}{{\mathbb D}}
\newcommand{\bF}{{\mathbb F}}
\newcommand{\bP}{{\mathbb P}}
\newcommand{\bQ}{{\mathbb Q}}

\newcommand{\boR}{{\mathbf R}}
\newcommand{\cF}{{\mathcal F}}

\newcommand{\cH}{{\mathcal H}}
\newcommand{\cX}{{\mathcal X}}
\newcommand{\oX}{\overline{X}}

\newcommand{\opi}{\overline{\pi}}
\newcommand{\ocX}{\overline{\mathcal X}}
\newcommand{\BM}{\text{\rm BM}}
\newcommand{\charsp}{\hbox{\rm char}\,}

\newcommand{\prim}{\text{\rm prim}}
\newcommand{\red}{\text{\rm red}}
\newcommand{\Gr}{\hbox{\rm Gr}}
\newcommand{\Hom}{\hbox{\rm Hom}}
\newcommand{\IC}{\hbox{\rm IC}}
\newcommand{\IH}{\hbox{\rm IH}}

\newcommand{\Ker}{\hbox{\rm Ker}}
\newcommand{\Imm}{\hbox{\rm Im}}
\newcommand{\Perv}{\hbox{\rm Perv}}
\newcommand{\Spec}{\hbox{\rm Spec}}
\newcommand{\Sing}{\hbox{\rm Sing}}

\begin{document}
\title[Weight filtration on the cohomology
of algebraic varieties]
{Weight filtration on the cohomology\\
of algebraic varieties}
\author{Masaki Hanamura}
\address{Mathematical Institute, Tohoku University, Sendai 980 JAPAN}
\email{hanamura@math.tohoku.ac.jp}
\author{Morihiko Saito}
\address{RIMS Kyoto University, Kyoto 606-8502 Japan}
\email{msaito@kurims.kyoto-u.ac.jp}
\date{\ver}
\begin{abstract}
We show that the etale cohomology (with compact supports) of an
algebraic variety $X$ over an algebraically closed field has
the canonical weight filtration $W$,
and prove that the middle weight part of the cohomology with
compact supports of $X$ is a subspace of the intersection cohomology
of a compactification $X'$ of X, or equivalently, the middle weight
part of the (so-called) Borel-Moore homology of $X$ is a quotient of
the intersection cohomology of $X'$.
We are informed that this has been shown by A. Weber in the case
$X$ is proper (and $k=\bC$)
using a theorem of G. Barthel, J.-P. Brasselet, K.-H. Fieseler,
O. Gabber and L. Kaup on morphisms between intersection
complexes.
We show that the assertion immediately follows from Gabber's purity
theorem for intersection complexes.
\end{abstract}
\maketitle

\bigskip
\centerline{\bf Introduction}

\bigskip\noindent
Let
$ X $ be an irreducible algebraic variety defined over an
algebraically closed field
$ k $.
Let
$ H_{j}^{\BM}(X) $ denote (so-called) Borel-Moore homology,
i.e. the dual of
$ H_{c}^{j}(X) $.
Here cohomology means etale cohomology with
$ l $-adic coefficients where
$ l $ is a prime number different from
$ \charsp k $ (and similarly for cohomology with
compact supports).
In the case
$ k = \bC $, we may also assume that it is singular cohomology with
$ \bQ $-coefficients as in [1]
using the comparison theorem.
If
$ \charsp k = 0 $ or
$ k $ is the algebraic closure of a finite field, it is known
that these cohomology and homology have the canonical weight 
filtration
$ W $ so that
$ H_{j}^{\BM}(X) $ has weights
$ \ge -j $,
see [2], [5], [7].
However, the general case of characteristic $p>0$ does not seem
to be treated in the references.
In this paper, we show that this easily follows from [2], [6]
using a model, see (1.2).

Let
$ \oX $ be a compactification of
$ X $, and
$ \IH^{j}(\oX) $ denote intersection cohomology, see [2], [10].
Note that
$ \IH^{j}(\oX) $ is pure of weight
$ j $ by Gabber's purity theorem [2] together with the stability of
pure complexes by the direct image under a proper morphism [7].
We have a canonical morphism
$$
\IH^{2n-j}(\oX)(n) \to H_{j}^{\BM}(X),
$$
which is the dual of
$ H_{c}^{j}(X) \to \IH^{j}(\oX) $, where
$ n = \dim X $, and
$ (n) $ is the Tate twist.
It was shown by G.~Barthel, J.-P.~Brasselet, K.-H.~Fieseler, O.~Gabber
and L.~Kaup (see [1]) that any algebraic cycle classes can be lifted
(noncanonically) to the intersection cohomology by the above morphism
at least when
$ k = \bC $.
In this paper we prove

\medskip\noindent
{\bf Theorem 1.} {\it
The canonical morphism
$ \IH^{2n-j}(\oX)(n) \to \Gr_{W}^{-j}H_{j}^{\BM}(X) $ is surjective, or
equivalently,
$ \Gr_{W}^{j}H_{c}^{j}(X) \to \IH^{j}(\oX) $ is injective.
}

\medskip
This is a refinement of the theorem in [1] explained above, since the
cycle class of an algebraic
$ d $-cycle is defined in
$ (\Gr_{W}^{-2d}H_{2d}^{\BM}(X))(-d) $.
We are informed that the assertion has been proved by A.~Weber [22]
in the case $X$ is proper (and $k=\bC$)
using a theorem in [1] on morphisms between intersection complexes.
In the case
$ \oX $ is smooth and $k=\bC$, Theorem~1 follows from the construction of
mixed Hodge structure in [5].
In general we have the surjectivity of the canonical morphism
$$
\IH^{2n-j}(\oX)(n)\to\Gr_{-j}^{W}(\IH^{2n-j}(X)(n)),
$$
factoring the first morphism in Theorem 1,
see Remark (i) in (2.2).
The proof of Theorem 1 immediately follows from Gabber's purity theorem
for intersection complexes ([2], [3]).
We also give a shorter proof of the theorem in [1] on morphisms
between intersection complexes (here it is not necessary to assume
$ k = \bC $):

\medskip\noindent
{\bf Theorem 2} (see [1]).
{\it Let
$ f : X \to Y $ be a morphism of equidimensional varieties over
an algebraically closed field
$ k $.
Then there is a noncanonical morphism
$$
f^{*}\IC_{Y}\bQ_{l} \to \IC_{X}\bQ_{l},
$$
whose composition with the canonical morphism
$ \bQ_{l} \to f^{*}\IC_{Y}\bQ_{l} $ coincides with the
canonical morphism
$ \bQ_{l} \to \IC_{X}\bQ_{l} $.
}

\medskip
In [1], this was proved by reducing to the case where
$ X $ is a closed subvariety of codimension
$ 1 $ in
$ Y $, and using induction on stratum, see also [21].
We show that Theorem~2 immediately follows from Gabber's purity
Theorem.
In the case of a closed embedding of irreducible varieties with relative
dimension
$ 1 $, they also showed the existence of a canonical morphism,
see [1].
We give another proof of it using the weight filtration on mixed perverse
sheaves [2], see (2.4).

If
$ X, Y $ are irreducible,
$ f $ is proper surjective, and
$ \boR f_{*}\IC_{X}\bQ_{l} $ is a shifted perverse sheaf on
$ Y $ (e.g. if
$ \dim X = \dim Y = 2 $ or
$ f $ is finite and surjective), then the semisimplicity of pure
perverse sheaves [2] implies the canonical morphisms
$$
\IC_{Y}\bQ_{l} \to \boR f_{*}\IC_{X}\bQ_{l},\quad
\boR f_{*}\IC_{X}\bQ_{l} \to \IC_{Y}\bQ_{l}.
$$
This can be used to construct cohomological correspondences, see [8].

This paper is a consequence of our discussions on the first author's
work [14] (see also [4]) at the conferences at Kagoshima (December 2004)
and at Baltimore (March 2005).
We thank the organizers of the conferences.
We also thank H.~Esnault and L.~Illusie for useful comments and questions
about the canonical choice of the morphisms, and
O.~Gabber for pointing out an error in an early version of this paper.

In Section 1 we show the existence of the weight filtration on the cohomology
of a variety over any algebraically closed field.
In Section 2 we prove the main theorems using Gabber's purity theorem.
In Section 3 we give an application of the decomposition theorem
of Beilinson, Bernstein and Deligne.

\newpage
\centerline{{\bf 1. Preliminaries}}

\bigskip\noindent
{\bf 1.1.~Intersection cohomology.}
Let
$ X $ be an algebraic variety (i.e. a separated reduced scheme of
finite type) over an algebraically closed field.
Let
$ l $ be a prime different from the characteristic of
$ k $.
By [2], [7], we have the bounded derived category with constructible
cohomology
$ D_{c}^{b}(X,\bQ_{l}) $.
Then the category of perverse sheaves
$ \Perv(X,\bQ_{l}) $ is defined to be a full subcategory of
$ D_{c}^{b}(X,\bQ_{l}) $ satisfying certain conditions, and
it is an abelian category, see [2].

If
$ X $ is purely
$ n $-dimensional, there is uniquely a shifted perverse sheaf
$ \IC_{X}\bQ_{l} \in D_{c}^{b}(X,\bQ_{l}) $ such that
$ (\IC_{X}\bQ_{l})[n] \in \Perv(X,\bQ_{l}) $, the restriction of
$ \IC_{X}\bQ_{l} $ to the smooth part
$ U $ of
$ X $ is the constant sheaf
$ \bQ_{l} $, and
$ (\IC_{X}\bQ_{l})[n] $ has no nontrivial sub nor quotient objects in
$ \Perv(X,\bQ_{l}) $ which are supported on the complement of
$ U $, see [2].
This is called the intersection complex.
To simplify the notation, we use the above normalization of
intersection complex as in [1], which is different from [2].
We have unique morphisms
$$
\bQ_{l} \to \IC_{X}\bQ_{l},\quad\IC_{X}\bQ_{l}(n)[2n] \to \bD\bQ_{l}
\quad\text{on}\quad X,
\leqno(1.1.1)
$$
which are dual of each other, and
whose restriction to the smooth part
$ U $ of
$ X $ is the identity on
$ \bQ_{l}|_{U} $ or
$ \bD\bQ_{l}|_{U} $.
Here
$ (n) $ denotes the Tate twist,
$ \bQ_{l} $ is a constant sheaf on
$ X $, and
$ \bD\bQ_{l} $ is the dual of
$ \bQ_{l} $ which is isomorphic to the dualizing complex
$ a^{!}\bQ_{l} $ where
$ a : X \to \Spec\, k $.
The first morphism of (1.1.1) follows from the fact that
$ {}^{p}\cH^{j}\bQ_{l} = 0 $ for
$ j > n $ and
$ (\IC_{X}\bQ_{l})[n] $ is a canonical quotient of
$ {}^{p}\cH^{n}\bQ_{l} $.
For the last assertion, note that
$ {}^{p}\cH^{n}\bQ_{l} $ has no nontrivial quotient
$ M $ whose support has dimension
$ < n $, where
$ {}^{p}\cH^{j} $ denotes the perverse cohomology sheaf.
(Indeed, this follows from
$$
\Hom({}^{p}\cH^{n}\bQ_{l},M) = \Hom(\bQ_{l}[n],M) =
\Hom(i^{*}\bQ_{l}[n],i^{*}M) = 0,
$$
by the vanishing of negative extension groups for perverse sheaves [2],
where
$ i $ denotes the inclusion of the support of
$ M $ into
$ X $.)

We define intersection cohomology (with
compact supports) by
$$
\IH^{j}(X,\bQ_{l}) = H^{j}(X,\IC_{X}\bQ_{l}),\quad
\IH_{c}^{j}(X,\bQ_{l}) = H_{c}^{j}(X,\IC_{X}\bQ_{l}).
$$
We will denote their dual by
$\IH_{j}(X,\bQ_{l})$ and $\IH^{\BM}_{j}(X,\bQ_{l})$ respectively.
These are isomorphic to
$\IH_c^{2n-j}(X,\bQ_{l})(n)$ and $\IH^{2n-j}(X,\bQ_{l})(n)$ respectively.

By (1.1.1) we have canonical morphisms
$$
H_{c}^{j}(X,\bQ_{l}) \to \IH_{c}^{j}(X,\bQ_{l}),\quad
\IH^{2n-j}(X,\bQ_{l})(n) \to H_{j}^{\BM}(X,\bQ_{l}),
\leqno(1.1.2)
$$
which are dual of each other.
Here
$ H_{j}^{\BM}(X,\bQ_{l}) := H^{-j}(X,\bD\bQ_{l}) $ denotes (so-called)
Borel-Moore homology which is the dual of
$ H_{c}^{j}(X,\bQ_{l}) $.

\medskip\noindent
{\bf 1.2.~Weight filtration.}
If
$ k $ is an algebraic closure of a finite field, then there is the
weight filtration
$ W $ on
$ H^{j}(X,\bQ_{l}) $,
$ H_{c}^{j}(X,\bQ_{l}) $ and hence on their dual
$ H_{j}(X,\bQ_{l}) $,
$ H_{j}^{\BM}(X,\bQ_{l}) $.
Moreover,
$ H_{c}^{j}(X,\bQ_{l}) $ and
$ H_{j}^{\BM}(X,\bQ_{l}) $ have weights
$ \le j $ and
$ \ge -j $ respectively, i.e.
$ \Gr_{i}^{W}H_{c}^{j}(X,\bQ_{l})=\Gr_{-i}^{W}H_{j}^{\BM}(X,\bQ_{l})=0 $
for
$ i > j $, see [2], [7].

It is known that these can be generalized to the case of any
algebraically closed field.
Indeed, if
$ k $ has characteristic
$ 0 $,
then this is well known as part of the theory of mixed
motives via realizations, see [17].
The argument is similar for
$ \charsp k > 0 $.
In this case, we may assume that
$ X $ is defined over a subfield
$ k' $ of
$ k $ which is finitely generated over
$ \overline{\bF}_{q} $ where
$ \overline{\bF}_{q} $ is the algebraic closure of a finite field
$ \bF_{q} $ in
$ k $.
Then there is a morphism
$ \pi : \cX \to S $ of algebraic varieties over
$ \overline{\bF}_{q} $ such that the function field of
$ S $ is
$ k' $ and the geometric generic fiber of
$ \pi $ over the geometric generic point of
$ S $ defined by the inclusion
$ k' \to k $ is identified with
$ X $.
Here we may assume that
$ R^{j}\pi_{*}\bQ_{l} $ and
$ R^{j}\pi_{!}\bQ_{l} $ are smooth sheaves
shrinking
$ S $ if necessary.

By [2] these smooth sheaves have the canonical weight filtration
$ W $, which is compatible with the weight filtration on
$ H^{j}(\cX_{s},\bQ_{l}) $ and
$ H_{c}^{j}(\cX_{s},\bQ_{l}) $ for any closed point
$ s $ of
$ S $ using the generic base change theorem [6].
Indeed, the stalk at $s$ of the smooth pure perverse sheaf $\cF$ of
weight $r$ on $S$ is pure of weight $r-d$ with $d=\dim S$, because
$i_s^!\cF=i_s^*\cF\otimes i_s^!\bQ_l$ for any smooth $\bQ_l$-sheaf $\cF$
where $i_s:\Spec\,k\to S$ is defined by $s$.
Restricting to the stalk at the geometric generic point,
this induces the weight filtration
$ W $ on
$ H^{j}(X,\bQ_{l}) $,
$ H_{c}^{j}(X,\bQ_{l}) $, and then on
$ H_{j}(X,\bQ_{l}) $,
$ H_{j}^{\BM}(X,\bQ_{l}) $ by duality.
The obtained filtration is independent of the choice of
$ k', S $.
This is functorial for morphisms of algebraic varieties in the usual
way,
because any morphism over $k$ has a model over $S$ replacing $S$ if
necessary.

\medskip\noindent
{\bf 1.3.~Remark.}
Let
$ \pi : \cX \to S $ be as above.
Then the pull-back of the intersection complex
$ \IC_{\cX}\bQ_{l} $ by the canonical morphism
$ X \to \cX $ is naturally isomorphic to
$ \IC_{X}\bQ_{l} $ up to a shift of complexes.
This can be shown by using the intermediate direct image [2]
together with the generic base change theorem [6].
As a corollary we get the purity of intersection cohomology in
the proper case using the above construction of weight filtration.

\bigskip\bigskip\centerline{{\bf 2. Proof of the main theorems}}

\bigskip\noindent
{\bf 2.1.~Proof of Theorem 1.}
Note that the canonical morphisms in Theorem~1 is
induced by (1.1.2) together with the restriction morphism
(and its dual).
We first reduce the assertion to the case where
$ k $ is an algebraic closure of a finite field
$ k_{0} $.
In the case
$ \charsp k = p > 0 $,
this follows from Remark (1.3) taking a closed point of
$ S $, because the canonical morphism
$$
\Gr_{W}^{j}H_{c}^{j}(X,\bQ_{l}) \to \IH^{j}(\oX,\bQ_{l})
$$
can be extended to a morphism of local systems
$$
\Gr_{W}^{j}R^{j}\pi_{!}\bQ_{l} \to
R^{j}\opi_{!}(\IC_{\ocX}\bQ_{l}),
$$
where
$ \opi : \ocX \to S $ is a compactification of
$ \pi $ whose geometric generic fiber is
$ \oX $.
If
$ \charsp k = 0 $,
the argument is similar using the mod
$ p $ reduction argument in [2].
So we may assume that
$ k $ is the algebraic closure of a finite field.

Let
$ K $ denote the mapping cone
$ C(j'_{!}\bQ_{l} \to \IC_{\oX}\bQ_{l}) $, where
$ j' : X \to \oX $ denotes the inclusion morphism.
By the associated long exact sequence
$$
\to H^{j-1}(\oX,K) \to H_{c}^{j}(X,\bQ_{l}) \to
\IH^{j}(\oX,\bQ_{l}) \to,
$$
it is sufficient to show that
$ H^{j-1}(\oX,K) $ has weights
$ \le j-1 $.
By [2], 5.1.14, this is reduced to the assertion that
$ K $ has weights
$ \le 0 $ in the sense of [2], 5.1.8.
By definition this is equivalent to the condition that
$ \cH^{j}K_{x} $ has weights
$ \le j $ for any closed point
$ x $ of
$ \oX $.
But this is further equivalent to the same assertion with
$ K $ replaced by
$ \IC_{\oX}\bQ_{l} $, because
$ \cH^{j}K_{x} = \cH^{j}(\IC_{\oX}\bQ_{l})_{x}/\bQ_{l} $ if
$ x \in X $,
$ j = 0 $, and
$ \cH^{j}K_{x} = \cH^{j}(\IC_{\oX}\bQ_{l})_{x} $ otherwise.
(The isomorphism for
$ j = 0 $ follows from the fact that the constant sheaf
$ \bQ_{l} $ has no nontrivial subsheaf supported on a closed
subvariety of strictly smaller dimension, which we apply to the
kernel of
$ \cH^{0} $ of the first morphism in (1.1.1).)
Then the above assertion for the intersection complex
$ \IC_{\oX}\bQ_{l} $ is known as
Gabber's purity theorem [2].
This completes the proof of Theorem~1.

\medskip\noindent
{\bf 2.2.~Remarks.}
(i) The canonical morphism
$ \IH^{j}(\oX)\to\Gr_{j}^{W}\IH^{j}(X) $ is surjective.
This follows from the distinguished triangle
$$
i_{*}i^{!}\IC_{\oX}\bQ_{l}\to\IC_{\oX}\bQ_{l}\to
\boR j'_{*}\IC_{X}\bQ_{l}\buildrel{+1}\over\to,
$$
(where
$ i : \oX\setminus X \to \oX $ is the inclusion of the complement)
together with the assertion that
$ i^{!}\IC_{\oX}\bQ_{l} $ has weights
$ \ge 0 $, see [2], 5.1.14.

\medskip
(ii) If
$ \oX $ is projective and
$ \charsp k = 0 $, then we have a canonical choice of a splitting of
the surjective or injective morphisms in Theorem~1 using a polarization
of Hodge structure.
However, it is not clear whether it is a good one.
For example, certain algebraic cycle classes (i.e. the graph of
`placid' maps) can be canonically lifted to the intersection cohomology
in [11] (see also [1]), and the relation with this is unclear.
This problem of canonical lifting is related to the Lefschetz
trace formula for intersection cohomology
(see [9], [11], [12], [13], [15], [16], [18], [19], [20], [23]),
and will be treated in [8].

\medskip\noindent
{\bf 2.3.~Proof of Theorem 2.}
Let
$ K = C(\bQ_{l} \to \IC_{Y}\bQ_{l}) $.
Using the distinguished triangle
$$
f^{*}K[-1] \to \bQ_{l} \to f^{*}\IC_{Y}\bQ_{l} \buildrel{+1}\over\to,
\leqno(2.3.1)
$$
the assertion is
equivalent to the vanishing of the composition of canonical morphisms
$$
f^{*}K[-1] \to \bQ_{l} \to \IC_{X}\bQ_{l}\quad \hbox{in}
\,\,\, D_{c}^{b}(X,\bQ_{l}).
\leqno(2.3.2)
$$
Then the assertion is reduced to the case where
$ k $ is an algebraic closure of a finite field.
Indeed, in the positive characteristic case, we have
$$
\Hom(f^{*}K[-1],\IC_{X}\bQ_{l}) =
H^{1}(X,\boR\cH{om}(f^{*}K,\IC_{X}\bQ_{l})),
$$
and this vector space can be extended to a
local system (i.e. a smooth sheaf) on
$ S $ as in (1.2).
Note that (2.3.2) is extended to a section of this local system and
its vanishing is equivalent to that for its restriction
over a closed point using the generic base change theorem [6].
(If
$ \charsp k > 0 $, it is also possible to replace
$ X $ with
$ \cX $ in (1.3).)
In the case
$ \charsp k = 0 $, this follows from the mod
$ p $ reduction argument in [2].

By the same argument as in (2.1),
$ f^{*}K $ has weights
$ \le 0 $.
There is a finite subfield
$ k_{0} $ of
$ k $ together with a morphism
$ f_{0} : X_{0} \to Y_{0} $ of algebraic varieties over
$ k_{0} $ such that
$ f $ is the base change of
$ f_{0} $ and the above morphisms are defined over
$ X_{0} $.
Then the assertion follows from [2], 5.1.15 (applied to
$ f^{*}K $ and
$ \IC_{X}\bQ_{l} $ which have weights
$ \le 0 $ and
$ \ge 0 $ respectively with the normalization of intersection
complexes in [1]).
This completes the proof of Theorem 2.

\medskip
The following was proved in [1] using induction on stratum at least if
$ \charsp k = 0 $.
We give another proof using the weight filtration
$ W $ on mixed perverse sheaves [2].

\medskip\noindent
{\bf 2.4.~Proposition.}
{\it Assume
$ X $,
$ Y $ irreducible,
$ f : X \to Y $ is a closed immersion, and
$ \dim Y = \dim X + 1 \,(= n + 1) $.
Then there is a canonical choice of the morphism
$f^{*}\IC_{Y}\bQ_{l} \to \IC_{X}\bQ_{l}$ satisfying the
condition in Theorem~2.
}

\medskip\noindent
{\it Proof.}
Let
$ j : Y \setminus X \to Y $ denote the inclusion morphism.
We have the long exact sequence
$$
\to {}^{p}\cH^{i}\IC_{Y}\bQ_{l}
\to f_{*}{}^{p}\cH^{i}f^{*}\IC_{Y}\bQ_{l}
\to {}^{p}\cH^{i+1}j_{!}j^{*}\IC_{Y}\bQ_{l}
\to,
$$
where
$ {}^{p}\cH^{i+1}j_{!}j^{*}\IC_{Y}\bQ_{l} = 0 $ for
$ i+1 > n+1 \,(= \dim Y) $, see [2], 4.2.4.
So we get
$$
{}^{p}\cH^{i}f^{*}\IC_{Y}\bQ_{l} = 0\quad\text{for}\,\,\,
i > n,
$$
because the perverse sheaf
$ \IC_{Y}\bQ_{l}[n+1] $ has no nontrivial quotient whose support
is contained in
$ X $.
Furthermore,
$ f^{*}\IC_{Y}\bQ_{l}[n] $ and hence
$ {}^{p}\cH^{n}f^{*}\IC_{Y}\bQ_{l} $ (see [2], 5.4.1) have weights
$ \le n $.
Thus we get a canonical morphism
$$
f^{*}\IC_{Y}\bQ_{l}[n]\to\Gr_{n}^{W}{}^{p}\cH^{n}f^{*}\IC_{Y}\bQ_{l},
$$
(factoring through
$ {}^{p}\cH^{n}f^{*}\IC_{Y}\bQ_{l} $).
By the semisimplicity theorem ([2], 5.3.8),
the target of the above morphism is semisimple, and is a direct sum of
a subobject whose support is strictly smaller than
$ X $ and an intersection complex with support
$ X $ associated with a smooth sheaf defined on a dense open subvariety of
$ X $.
So a morphism
$$
f^{*}\IC_{Y}\bQ_{l} \to \IC_{X}\bQ_{l},
\leqno(2.4.1)
$$
is uniquely determined by its restriction to any sufficiently small
nonempty open subvariety
$ X' $ of
$ X $ (contained in the smooth part of
$ X $), because
$ \IC_{X}\bQ_{l}[n] $ is a perverse sheaf and is pure of weight
$ n $.
We have a similar assertion for a morphism
$ \bQ_{l} \to \IC_{X}\bQ_{l} $.

Let
$ f' ; X' \to Y' $ be the restriction (or base change) of
$ f $ over an open subvariety
$ Y' $ of
$ Y $ whose complement has codimension
$ \ge 2 $.
Let
$ \rho' : \widetilde{Y}' \to Y' $ be the normalization, and
$ X'' = (\widetilde{Y}' \times_{Y'}X')_{\red} $ with the canonical
morphism
$ \rho'' : X'' \to X ' $.
Let
$ X''_{i} $ be the irreducible components of
$ X'' $ with
$ s_{i} $ the separable degree of the extension
$ k(X''_{i})/k(X') $.
Here we may assume that
$ \widetilde{Y}' $ and
$ X'' $ are smooth over
$ k $, and
$ f'{}^{*}\rho'_{*}\bQ_{l} = \rho''_{*}\bQ_{l} $ is a local system
(shrinking
$ Y' $ if necessary).
Then we have
$$
\IC_{Y'}\bQ_{l} = \rho'_{*}\bQ_{l},\quad
f^{*}\IC_{Y'}\bQ_{l} = \rho''_{*}\bQ_{l},
$$
and the restriction of the morphism (2.4.1) to
$ X' $ is given by
$ \rho''_{*}\bQ_{l} \to \bQ_{l} $.

Let
$ X_{i}^{s} $ be the irreducible variety with finite morphisms
$ X''_{i} \to X_{i}^{s} \to X' $ factoring
$ \rho'' $ so that the function field
$ k(X_{i}^{s}) $ is the maximal separable extension of
$ k(X') $ contained in
$ k(X''_{i}) $ (shrinking
$ X' $ if necessary).
Then the direct image of the constant sheaf
$ \bQ_{l} $ by
$ X''_{i} \to X_{i}^{s} $ is the constant sheaf
$ \bQ_{l} $ so that we may replace
$ \rho''_{*}\bQ_{l}$ with
$ \rho^{s}_{*}\bQ_{l}$ where
$ \rho^{s} $ is the canonical morphism of the disjoint union of the
$ X_{i}^{s} $ to
$ X' $.
So the desired morphism
$ \rho^{s}_{*}\bQ_{l} \to \bQ_{l} $ is given by the dual of the
canonical morphism
$ \iota : \bQ_{l} \to \rho^{s}_{*}\bQ_{l} $, divided by
$ \sum_{i}s_{i} $, because the composition of
$ \iota $ and its dual
$ \iota^{\vee} : \rho^{s}_{*}\bQ_{l} \to \bQ_{l} $ is the
multiplication by
$ \sum_{i}s_{i} $ on
$ \bQ_{l} $.
This completes the proof of Proposition~(2.4).

\medskip\noindent
{\bf 2.5.~Remarks.} (i)
Even in the case of Proposition (2.4), the morphism satisfying the
condition in Theorem~2 is not unique in general (for example, if
$ Y $ has etale locally two components whose intersection is
$ X $).
However, it is unique if
$ Y $ is etale locally irreducible at the generic point of
$ X $, see also [1].

\medskip
(ii) In the case
$ f $ is a closed immersion of codimension
$ \ge 2 $, the morphism satisfying the condition in Theorem~2 is not
unique even if
$ X $,
$ Y $ have only isolated singularities.
For simplicity, assume
$ \dim Y = \dim X + 2 \,(= n + 2) $.
For
$ x \in \Sing \,X $, let
$$
E_{X,x}^{i} = \cH^{i}(\IC_{X}\bQ_{l})_{x},\quad
E_{Y,x}^{i} = \cH^{i}(\IC_{Y}\bQ_{l})_{x}.
$$
Then we have a distinguished triangle
$$
{}^{p}\tau_{\le n}f^{*}\IC_{Y}\bQ_{l} \to f^{*}\IC_{Y}\bQ_{l} \to
E_{Y,x}^{n+1}[-n-1]\buildrel{+1}\over\longrightarrow,
$$
and there is a contribution of
$$
\Hom(E_{Y,x}^{n+1}[-n-1],\IC_{X}\bQ_{l}),
$$
to the
ambiguity of the morphism using the above distinguished triangle
together with the vanishing of negative extensions.
Moreover, by the self duality
$$
\bD \IC_{X}\bQ_{l} = \IC_{X}\bQ_{l}(n)[2n],
$$
the last group is isomorphic to
$$
\aligned
&\Hom(\IC_{X}\bQ_{l},(\bD E_{Y,x}^{n+1})(-n)[-n+1])
\\
&\quad =\Hom(E_{X,x}^{n-1},(\bD E_{Y,x}^{n+1})(-n)),
\endaligned
$$
which is not necessarily zero in general.
For example, if
$ X $ is the affine 
cone of a smooth projective variety
$ V $, then
$ E_{X,x}^{j} $ is isomorphic to the primitive cohomology
$ H^{j}_{\prim}(V,\bQ_{l}) $.

\bigskip\bigskip
\centerline{\bf 3. Application of the decomposition theorem}

\bigskip\noindent
{\bf 3.1.~Decomposition theorem.}
Let $f:X\to Y$ be a proper surjective morphism of irreducible algebraic
varieties over an algebraically closed field $k$. Let $n=\dim X$.
By the decomposition theorem [2], we have a noncanonical isomorphism
$$\boR f_*\IC_X\bQ_l[n]\simeq\mopls_{i,Z}\IC_Z E_{Z^{\circ}}^i
[\dim Z][-i]\quad\text{in}\quad D^b_c(Y,\bQ_l),\leqno(3.1.1)$$
where $Z$ are irreducible reduced closed subvarieties of $Y$, and
$E^i_{Z^{\circ}}$ are smooth $l$-adic sheaves defined on dense open
smooth subvarieties
$Z^{\circ}$ of $Z$ which we may assume to be independent of $i$.
This can be reduced to the case where the base field
$ k $ is finite by using a model as in (1.2) together with [6].
Note that the shift by $n$ or $\dim Z$ is needed for 
$\IC_X\bQ_l$, $\IC_Z E_{Z^{\circ}}^i$ because of the normalization
of the intersection complex in this paper.
Set $$m_Z=\max\{i\,|\,E_{Z^{\circ}}^i\ne 0\}=
\max\{i\,|\,E_{Z^{\circ}}^{-i}\ne 0\}.$$
Note that $E_{Z^{\circ}}^i$ is the dual of $E_{Z^{\circ}}^{-i}$
up to a Tate twist.

\medskip\noindent
{\bf 3.2.~Lemma.} {\it
With the above notation, assume $X$ is smooth over $k$ so that
$\IC_X\bQ_l=\bQ_l$. Then
}
$$m_Z\le n-\dim Z-2\,\,\,\hbox{if}\,\,\,Z\ne Y,\quad m_Y\le n-\dim Y.
\leqno(3.2.1)$$

\medskip\noindent
{\it Proof.}
Let $X_s=f^{-1}(s)$.
For a general closed point $s$ of $Z^{\circ}$, we have
$$m_Z\le 2\dim X_s+\dim Z-n,$$
because the stalk of $R^jf_*(\bQ_l[n])$ at $s$
vanishes unless $j\in[-n,2\dim X_s-n]$ and the stalk of
$\IC_Z E_{Z^{\circ}}^i[\dim Z][-i]$ at $s$ is $E^i_{Z^{\circ},s}$ put at
the degree $j=i-\dim Z$.
We have $\dim X_s\le n-\dim Z-1$ if $Z\ne Y$. So the assertion follows.

\medskip
As a corollary we have the following
(which is used in [8]).

\medskip\noindent
{\bf 3.3.~Proposition.} {\it
With the notation of $(3.1)$, assume $X$ smooth.
Let $$u:\bQ_l[n]\to \boR f_*\bQ_l[n],\quad
v:\boR f_*\bQ_l[n]\to\bD\bQ_l(-n)[-n]$$ be the canonical morphisms on $Y$,
where $v$ is the dual of $u$. Let
$$\aligned u_{Z,i}&:\bQ_l[n]\to\IC_Z E_{Z^{\circ}}^i[\dim Z][-i],\\
v_{Z,i}&:\IC_Z E_{Z^{\circ}}^i[\dim Z][-i]\to\bD\bQ_l(-n)[-n]\endaligned
\leqno(3.3.1)$$
be the induced morphisms using the decomposition $(3.1.1)$.
Then
}
$$\aligned u_{Z,i}=0\,\,\,\text{for}\,\,\,(Z,i)\ne(Y,n-\dim Y),\\
v_{Z,i}=0\,\,\,\text{for}\,\,\,(Z,i)\ne(Y,\dim Y-n).\endaligned
\leqno(3.3.2)$$

\medskip\noindent
{\it Proof.}
By the adjunction for the inclusion $Z\to Y$, we have the canonical
isomorphism
$$\aligned\Hom_Y&(\bQ_l[n],\IC_Z E_{Z^{\circ}}^i[\dim Z][-i])\\
=\Hom_Z&(\bQ_l[\dim Z],\IC_Z E_{Z^{\circ}}^i[\dim Z][\dim Z-n-i]),
\endaligned$$
and the last group vanishes for $\dim Z-n-i<0$
by the semi-perversity of the constant sheaf $\bQ_l[\dim Z]$ on $Z$.
So we get the first assertion by (3.2.1).
The assertion is similar for the second.

\medskip
As a corollary of (3.3), we have

\medskip\noindent
{\bf 3.4.~Corollary.} {\it
Let $f:X\to Y$ be a proper surjective morphism of irreducible algebraic
varieties over an algebraically closed field $k$.
Assume $X$ is smooth over $k$.
Then
}
$$\aligned\Ker(H^i(Y,\bQ_l)\to H^i(X,\bQ_l))&=
\Ker(H^i(Y,\bQ_l)\to \IH^i(Y,\bQ_l)),\\
\Ker(H_c^i(Y,\bQ_l)\to H_c^i(X,\bQ_l))&=
\Ker(H_c^i(Y,\bQ_l)\to\IH_c^i(Y,\bQ_l)),\\
\Imm(H_i(X,\bQ_l)\to H_i(Y,\bQ_l))&=
\Imm(\IH_i(Y,\bQ_l)\to H_i(Y,\bQ_l)),\\
\Imm(H^{\BM}_i(X,\bQ_l)\to H^{\BM}_i(Y,\bQ_l))&=
\Imm(\IH^{\BM}_i(Y,\bQ_l)\to H^{\BM}_i(Y,\bQ_l)).\endaligned$$

\medskip\noindent
{\it Proof.}
This follows from (3.3).

\medskip\noindent
{\bf 3.5.~Remarks.}
(i) The above corollary is related to a question of A.~Weber
when $X\to Y$ is a resolution of singularities.

\medskip
(ii) The first isomorphism of (3.4) implies a proof of Theorem~1 in the
case $k=\bC$ and $X$ proper, using the weight spectral sequence in [5].
Indeed, if $X\to Y$ is a resolution of singularities, then
$\Gr^W_iH^i(Y)$ is a subspace of $H^i(X)$ (which coincides with the image
of $H^i(Y)$) by the weight spectral sequence associated to a simplicial
resolution as in Remark (iii) below.
So we can replace $H^i(X)$ with $\IH^i(Y)$ by (3.4).

\medskip
(iii)
Let $Y$ be a proper irreducible variety of dimension $n$ over an
algebraically closed field $k$.
Consider the natural morphisms
$$\Gr_n^W H^n(Y,\bQ_l)\buildrel{u_n}\over\longrightarrow
\IH^n(Y,\bQ_l)\buildrel{v_n}\over\longrightarrow
\Gr_n^W(H_n(Y,\bQ_l)(-n)),$$
where $u_n$ is injective and $v_n$ is the dual of $u_n$
and is surjective.
These are induced by the morphisms
$u,v$ in (3.3).
We have a canonical self-pairing of $\IH^n(Y,\bQ_l)$,
and a canonical pairing between $H^n(Y,\bQ_l)$ and
$H_n(Y,\bQ_l)(-n)$ with values in $\bQ_l(-n)$ so that
$$ \langle u(\xi),\eta\rangle=\langle \xi,v(\eta)\rangle\quad
\text{for}\,\,\,\xi\in H^n(Y,\bQ_l),v\in\IH^n(Y,\bQ_l).\leqno(3.5.1)$$
By (3.3) and (3.4) we can replace
$\IH^n(Y,\bQ_l)$ with $H^n(X,\bQ_l)$
where $X$ is a resolution of singularities of $Y$.
By (3.1.1) $\IH^n(Y,\bQ_l)$ is a direct factor of $H^n(X,\bQ_l)$
noncanonically,
and the morphisms $u$ and $v$ are compatible with
any decomposition (3.1.1) by (3.3).

\medskip
(iv) With the above notation, the restriction of the self-pairing to
the image of $u_n$ does not seem to be nondegenerate in general.
Indeed, if we take a simplicial resolution $X_{\ssbull}$
of $Y$ such that $X_0=X$ in the case of $k=\bC$,
then this image coincides with the kernel of $H^n(X)\to H^n(X_1)$.
However, this does not seem to be compatible with the Lefschetz
decomposition in general.

For example, let $X$ be the blow-up of $\bP^2$ along a point.
This is a $\bP^1$-bundle over $C=\bP^1$ having disjoint two sections
$C_i$ whose self-intersection number is $i$ for $i=\pm 1$.
If there is a variety $Y$ which is obtained by identifying the
two sections $C_1$ and $C_{-1}$, then the restriction of the canonical
pairing to the kernel of $H^2(X)\to H^2(C)$ would be degenerate,
where the last morphism is defined by the difference between the
restriction morphisms to $C_1$ and $C_{-1}$ (both identified with $C$
naturally).
However, it is not clear if such $Y$ exists in the category of algebraic
varieties (although it exists as an analytic space if $k=\bC$).
It would not be projective at least.

\medskip\noindent
{\bf 3.6.~Remark.}
Let
$ f : X \to Y $ be a proper surjective morphism of irreducible
varieties over an algebraically closed field $k$ as in (3.1).
Let $n=\dim X$.
Assume
$$
\boR f_{*}\IC_{X}\bQ_{l}[n] \,\,\text{is a perverse sheaf on}
\,\, Y.
\leqno(3.6.1)
$$
This assumption implies that $\dim X=\dim Y=n$, and the direct sum
decomposition (3.1.1) becomes the canonical decomposition in the category
of perverse sheaves
$$
\boR f_{*}\IC_{X}\bQ_{l}[n] = \mopls_{Z}\IC_{Z}E_{Z^{o}}
[\dim Z],\leqno(3.6.2)
$$
because there are no nontrivial morphisms between intersection
complexes with different supports.
This induces canonical morphisms
$$
\IC_{Y}\bQ_{l} \to \boR f_{*}\IC_{X}\bQ_{l},\quad
\boR f_{*}\IC_{X}\bQ_{l} \to \IC_{Y}\bQ_{l}.
$$
Indeed, over a sufficiently small non-empty open subvariety of
$ Y $, the intersection complexes coincide with the constant
sheaf
$ \bQ_{l} $, and the assertion is clear.
Then we can extend the obtained morphisms uniquely over
$ Y $ using the decomposition (3.6.2) together with the intermediate
direct image [2].
Here we can neglect
$ \IC_{Z}E_{Z^{o}} $ for
$ Z \ne Y $, because an intersection complex has no nontrivial sub nor
quotient objects with strictly smaller support.

The condition (3.6.1) is satisfied for example if
$ \dim X = \dim Y = 2 $, or
$ f $ is finite and surjective.
In the first case, (3.6.1) follows from the fact that
the support of
$ \cH^{1}\IC_{X}\bQ_{l} $ is discrete.
In the second case, the direct image is an intersection complex,
because the direct image by a finite morphism is an exact functor
of perverse sheaves, see [2].

\end{document}